\documentstyle[amssymb,12pt]{article}

\begin{document}

\newtheorem{theorem}{Theorem}{}
\newtheorem{lemma}[theorem]{Lemma}{}
\newtheorem{corollary}[theorem]{Corollary}{}
\newtheorem{conjecture}[theorem]{Conjecture}{}
\newtheorem{proposition}[theorem]{Proposition}{}
\newtheorem{axiom}{Axiom}{}
\newtheorem{remark}{Remark}{}
\newtheorem{example}{Example}{}
\newtheorem{exercise}{Exercise}{}
\newtheorem{definition}{Definition}{}

\title{Regular nilpotent elements and quantum groups}
\author{
Alexey Sevostyanov \footnote{e-mail seva@teorfys.uu.se}\\ 
Institute of Theoretical Physics, Uppsala University}

\maketitle
\begin{flushright}
UU--ITP 4/98

\end{flushright}

\begin{abstract}
We suggest new realizations of quantum groups $U_q({\frak g})$
 corresponding to
complex simple Lie algebras, and of affine quantum groups.
These new realizations are labeled by Coxeter elements 
of the corresponding Weyl group and have the following
key feature: The natural counterparts of the subalgebras
$U({\frak n})$, where ${\frak n} \subset {\frak g}$ is
a maximal nilpotent subalgebra, possess non--singular
characters. 
\end{abstract}

\section*{Introduction}
Let $\frak g$ be a complex simple Lie algebra, 
$\frak b$ a Borel subalgebra, and ${\frak n}=[{\frak b},{\frak b}]$
its nilradical. We denote by  $(~,~)$ the Killing form on $\frak g$.
An element $f\in {\frak n}$ is called {\em regular nilpotent}
if its centralizer in $\frak g$ is of minimal possible dimension.
Any regular nilpotent element $f\in {\frak n}$ defines a character
$\chi$ of the opposite nilpotent subalgebra 
$\overline{\frak n}=[\overline{\frak b},\overline{\frak b}]$,
where $\overline{\frak b}$ is the opposite Borel subalgebra.
Naturally, the character $\chi$ extends to a character of 
the universal enveloping algebra $U(\overline{\frak n})$. 
Recall that $U(\overline{\frak n})$ is generated  by the 
positive simple root generators $X_i^+, i=l,\ldots ,rank~{\frak g}$
of the Chevalley basis associated with the pair 
$({\frak g},\overline{\frak b})$. On  these generators
the character $\chi$ takes values $c_i \neq 0$, 
$\chi (X_i^+)=c_i$.
Conversely, any  character of this form determines a regular
nilpotent element in ${\frak n}$. 

Regular nilpotent elements are of great importance  
in the structure theory of Lie algebras and in its applications.
In particular, there is a relation between regular nilpotent
elements of a complex semisimple Lie algebra
and Coxeter elements of the corresponding Weyl group \cite{K}.
Other applications of regular nilpotent elements
include  the theory of Whittacker modules in representation theory 
of semisimple Lie algebras \cite{K1}, the integrability
of the Toda lattice \cite{K2}, and the  remarkable realization
of the center of the universal enveloping algebra of a complex
simple Lie algebra as a Hecke algebra \cite{K1}.
This provides a motivation to look for counterparts
of regular nilpotent elements in the theory of quantum groups.

Let $U_q({\frak g})$ be the quantum group associated with a 
complex simple Lie algebra ${\frak g}$, and let $U_q({\frak n})$
be the subalgebra of $U_q({\frak g})$ corresponding to the nilpotent
Lie subalgebra ${\frak n} \subset {\frak g}$. $U_q({\frak n})$ is
generated by simple positive root generators of $U_q({\frak g})$
subject to the q--Serre relations. It is easy to show that $U_q({\frak n})$
has no nondegenerate characters (taking nonvanishing values
on all simple root generators)! Our first main result
is the  family of new realizations of the
quantum group $U_q({\frak g})$, one for each Coxeter element
in the corresponding Weyl group. The counterparts of $U({\frak n})$,
which naturally arise in these new realizations of $U_q({\frak g})$,
do have non--singular characters. Thus, we get proper quantum
counterparts of $U({\frak n})$ and of its non--singular characters.
As a byproduct, we derive an interesting formula for the Caley
transform of a Coxeter element.

Next, we generalize our consideration to the case of affine
Lie algebras. Similar to the finite-dimensional situation,
the  subalgebra  $U_q({\frak n}((z)))$ in the affine quantum group
$U_q(\widehat{\frak g})$ naturally corresponding to ${\frak n}((z))$ has no characters taking nonvanishing values
on the quantum counterparts of the  loop generators of ${\frak n}((z))$.
Again, we introduce new realizations of $U_q(\widehat{\frak g})$,
labeled by Coxeter elements, such that the natural counterparts
of $U({\frak n}((z)))$ acquire such characters. 
Our realizations are variations of the Drinfeld's `new realization'
of affine quantum groups  \cite{nr}.
The perspective application of our construction is
the Drinfeld-Sokolov reduction for affine quantum groups.

\subsubsection*{Acknowledgements.}

The author would like to thank B.Enriquez , G. Felder ,  for useful 
discussions. I am also grateful to A. Alekseev and to M. Golenishcheva-Kutuzova 
for careful reading of the text.

\section{Non--singular characters and finite--dimensional quantum groups}

In this section we construct quantum counterparts of the principal nilpotent 
Lie subalgebras of complex simple Lie algebras and of their non--singular 
characters.

We follow the notation of \cite{Kac}.
Let ${\frak h}^*$ be an $l$--dimensional complex vector space, 
$a_{ij} , i,j=1,\ldots ,l$  a Cartan matrix of finite type , 
$\Delta \in {\frak h}^*$  the corresponding root system, and
$\{\alpha _1,...,\alpha _l\}$ the set of simple roots. 
Denote by $W$ the Weyl group of the root system $\Delta$, and
by $s_1,...,s_l \in W$  reflections corresponding to simple roots. 
Let $d_1,\ldots , d_l$ be coprime positive integers 
such that the matrix $b_{ij}=d_ia_{ij}$ is symmetric. 
There exists a unique non--degenerate $W$--invariant 
scalar product $\left( ,\right) $ on ${\frak h}^*$ 
such that $(\alpha_i , \alpha_j)=b_{ij}$. 

Let $\frak{g}$ be the complex simple Lie algebra associated to 
the Cartan matrix $a_{ij}$. Denote by $\frak{n}\subset {\frak{g}}$
the principal nilpotent subalgebra generated by the simple positive
root generators of the Chevalley basis. 
\begin{definition}
A character $\chi: {\frak n} \rightarrow {\Bbb C}$ is called non--singular
if and only if it takes non-vanishing values on all simple root
generators of ${\frak n}$.
\end{definition}
Note that any non--singular character is equivalent (up to a Lie
algebra automorphism of ${\frak n}$) to $\chi_0$ which takes value $1$
on each simple root generator. Any character of ${\frak n}$ naturally
extends to a character of the associative algebra $U({\frak n})$.
It is our goal to construct quantum counterparts of the algebra
$U({\frak n})$ and of the non--singular character $\chi_0$.

Let $q$ be a complex number, $0<|q|<1$. Put $q_i=q^{d_i}$. We 
consider the simply--connected rational form $U_q^R({\frak g})$ 
of the quantum group $U_q({\frak g})$ \cite{ChP}, Section 9.1.
This is an
associative algebra over $\Bbb C$ with
generators $X_i^\pm , L_i , L_i^{-1} , i=1,\ldots , l$ 
subject to the relations:

\begin{equation}\label{qgr}
\begin{array}{l}
L_iL_j=L_jL_i~ ,~ L_iL_i^{-1}=L_i^{-1}L_i=1 ,\\
 \\
L_iX_j^\pm L_i^{-1}=q_i^{\pm \delta _{i,j}}X_j^\pm , \\
 \\
X_i^+X_j^- -X_j^-X_i^+ = \delta _{i,j}{K_i -K_i^{-1} \over q_i -q_i^{-1}} , \\
 \\
K_i=\prod_{j=1}^lL_j^{a_{ji}} , \\
 \\
\mbox{and the q--Serre relations:} \\
\\
\sum_{r=0}^{1-a_{ij}}(-1)^r 
\left[ \begin{array}{c} 1-a_{ij} \\ r \end{array} \right]_{q_i} 
(X_i^\pm )^{1-a_{ij}-r}X_j^\pm(X_i^\pm)^r =0 ,~ i \neq j ,\\ \\
\mbox{ where }\\
 \\
\left[ \begin{array}{c} m \\ n \end{array} \right]_q={[m]_q! \over [n]_q![n-m]_q!} ,~ 
[n]_q!=[n]_q\ldots [1]_q ,~ [n]_q={q^n - q^{-n} \over q-q^{-1} }.
\end{array}
\end{equation}
The elements $X_i^\pm$ correspond to the simple positive (negative) 
root generators. We would like to show that the algebra spanned by 
$X_i^+ , i=1, \ldots , l$ does not admit characters
which take nonvanishing values on all generators $X_i^+$,
except for the case of $U_q(sl(2))$ when q-Serre relations
do not appear.
 
Suppose, $\chi$ is such a character, and $\chi(X_i)=c_i$. The
q-Serre relations are homogeneous and, hence, one can put
$c_i=1$ for all $i$ without loss of generality.  
By applying the character $\chi$ to the q-Serre relations
one obtains a family of identities,
\begin{equation} \label{false}
\sum_{r=0}^{1-a_{ij}}(-1)^r 
\left[ \begin{array}{c} 1-a_{ij} \\ r \end{array} \right]_{q_i} =0 , 
\, i \neq j.
\end{equation}
We claim that some of these relations fail for  the quantized universal enveloping algebra $U_q^R({\frak g})$ 
of any simple Lie algebra $\frak g$ , with the exception of ${\frak g}=sl(2)$. In a more general setting , relations (\ref{false}) 
are analysed in the following lemma.

\begin{lemma}\label{qbinom} 
The only rational solutions of equation
\begin{equation}\label{c1}
\sum_{k=0}^{m}(-1)^k 
\left[ \begin{array}{c} m \\ k \end{array} \right]_{t}
t^{kc}=0 , 
\end{equation}
where $t$ is a complex number, $0< |t| < 1$, are of the form
\begin{equation}\label{c2}
c=-m+1,-m+2, \ldots ,m-2,m-1.
\end{equation}
\end{lemma}
{\em Proof.}
According to the q--binomial theorem \cite{GR}, 

\begin{equation}\label{z}
\sum_{k=0}^{m}(-z)^k 
\left[ \begin{array}{c} m \\ k \end{array} \right]_{t} 
=\prod_{p=0}^{m-1}(1-t^{m-1-2p}z).
\end{equation}

Put $z=t^c$ in this relation. Then the l.h.s of (\ref{z}) coincides with 
the l.h.s. of (\ref{c1}). 
Now (\ref{z}) implies that $c=m-1-2p , p=0, \ldots ,m-1$
are the only rational solutions of (\ref{c1}).

Now we return to identities (\ref{false}).
Any Cartan matrix contains at least one off-diagonal element
equal to $-1$. Then, $m= 1- a_{ij} = 2$ and 
$c=\pm 1$, and lemma \ref{qbinom} implies
that some of identities  (\ref{false}) are false for
any simple Lie algebra, except for $sl(2)$. Hence , subalgebras of $U_q^R({\frak g})$ generated by 
$X_i^+$ do not possess non--singular characters.

It is our goal to construct subalgebras of $U_q^R({\frak g})$ which resemble the subalgebra $U({\frak n}) 
\subset U({\frak g})$ and possess non--singular
characters. 
Denote by $S_l$ the symmetric group of $l$ elements.
To any element $\pi \in S_l$ we associate a Coxeter element $s_{\pi}$ by the formula
$s_\pi =s_{\pi (1)}\ldots s_{\pi (l)}$.
For each Coxeter element $s_\pi$ we define an associative algebra $F_q^\pi $ generated by elements $e_i ,~ i=1, \ldots l$ subject to the relations :

\begin{equation}\label{fqpi}
\sum_{r=0}^{1-a_{ij}}(-1)^r q^{r c_{ij}^{\pi}}
\left[ \begin{array}{c} 1-a_{ij} \\ r \end{array} \right]_{q_i} 
(e_i )^{1-a_{ij}-r}e_j (e_i)^r =0 ,~ i \neq j , 
\end{equation}
where $c_{ij}^{\pi}=\left( {1+s_\pi \over 1-s_\pi }\alpha_i , \alpha_j \right)$ are matrix elements of the Caley transform of $s_\pi$ 
in the basis of simple roots.

\begin{proposition}\label{charf}
The map $\chi_q^\pi:F_q^\pi \rightarrow {\Bbb C} $ defined on the generators by 
$\chi_q^\pi(e_i)=1$ is a character of the algebra $F_q^\pi $.
\end{proposition}

To show that  $\chi_q^\pi$ is a character of $F_q^\pi $ it is sufficient to check that the  defining 
relations (\ref{fqpi}) belong to the kernel of $\chi_q^\pi$ ,i.e.
\begin{equation}\label{chifqpi}
\sum_{r=0}^{1-a_{ij}}(-1)^r q^{r c_{ij}^{\pi}}
\left[ \begin{array}{c} 1-a_{ij} \\ r \end{array} \right]_{q_i}=0 ,~ i \neq j . 
\end{equation}

As a preparation for the proof of proposition \ref{charf} we study the matrix elements of the Caley transform of 
$s_\pi$ which enter the definition of $F_q^\pi $.

\begin{lemma}\label{tmatrel}
The matrix elements of  ${1+s_\pi \over 1-s_\pi }$ are of the form :

\begin{equation}\label{matrel}
\left( {1+s_\pi \over 1-s_\pi }\alpha_i , \alpha_j \right)=
\varepsilon_{ij}^\pi b_{ij},
\end{equation}
where
\begin{equation}
\varepsilon_{ij}^\pi =\left\{ \begin{array}{ll}
-1 & \pi^{-1}(i) <\pi^{-1}(j) \\
0 & i=j \\
1 & \pi^{-1}(i) >\pi^{-1}(j) 
\end{array}
\right .
\end{equation}

\end{lemma}

{\em Proof.} (compare \cite{Bur} , Ch. V , \S 6 , Ex. 3).

First we calculate the matrix of the Coxeter element $s_\pi$ with respect to the 
basis of simple roots. We obtain this matrix in the form of the Gauss 
decomposition of the operator $s_\pi$.

Let $z_{\pi (i)}=s_\pi \alpha_{\pi (i)}$. Recall that $s_i(\alpha_j)=\alpha_j-a_{ji}\alpha_i$. 
Using this definition the elements $z_{\pi (i)}$ may be represented as:

$$
z_{\pi (i)}=y_{\pi (i)} -\sum_{k \geq i} a_{\pi (k) \pi (i)}y_{\pi (k)},
$$
where

\begin{equation}\label{y}
y_{\pi (i)}=s_{\pi (1)}\ldots s_{\pi (i-1)}\alpha_{\pi (i)}.
\end{equation}

Using matrix notation we can rewrite the last formula as follows:

\begin{equation}\label{2*}
\begin{array}{l}
z_{\pi (i)}=
(I+V)_{\pi (k) \pi (i)}y_{\pi (k)} , \\  \\ \mbox{ where } V_{\pi (k) \pi (i)}=
\left\{ \begin{array}{ll}
a_{\pi (k) \pi (i)} & k\geq i \\
0 & k < i
\end{array}
\right  .
\end{array}
\end{equation}

To calculate the matrix of the operator $s_\pi$ with respect to the basis of simple roots we have to express 
the elements $y_{\pi (i)}$ via the simple roots.
Applying the definition of  simple reflections to (\ref{y}) we can pull out the element $\alpha_{\pi (i)}$ to the right:

\[
y_{\pi (i)}=\alpha_{\pi (i)}-\sum_{k<i}a_{\pi (k) \pi (i)}y_{\pi (k)}.
\]

Therefore

\[
\alpha_{\pi (i)}=(I+U)_{\pi (k) \pi (i)}y_{\pi (k)} ~, \mbox{ where } U_{\pi (k) \pi (i)}=
\left\{ \begin{array}{ll}
a_{\pi (k) \pi (i)} & k<i \\
0 & k \geq i
\end{array}
\right .
\]

Thus

\begin{equation}\label{1*}
y_{\pi (k)}=(I+U)^{-1}_{\pi (j) \pi (k)}\alpha_{\pi (j)}.
\end{equation}

Summarizing (\ref{1*}) and (\ref{2*}) we obtain:

\begin{equation}
s_\pi \alpha_i=\left( (I+U)^{-1}(I-V) \right)_{ki}\alpha_k .
\end{equation}

This implies:

\begin{equation}\label{3*}
{1+s_\pi \over 1-s_\pi}\alpha_i=\left( {2I+U-V \over U+V}\right)_{ki}\alpha_k .
\end{equation}

Observe that $(U+V)_{ki}=a_{ki}$ and $(2I+U-V)_{ij}=-a_{ij}\varepsilon_{ij}^\pi$. 
Substituting these expressions into (\ref{3*}) we get :

\begin{eqnarray}
\left( {1+s_\pi \over 1-s_\pi }\alpha_i , \alpha_j \right) = 
-(a^{-1})_{kp}\varepsilon_{pi}^\pi a_{pi}b_{jk}=\\
-d_ja_{jk}(a^{-1})_{kp}\varepsilon_{pi}^\pi a_{pi}  = 
\varepsilon_{ij}^\pi b_{ij}.
\end{eqnarray}

This concludes the proof of the lemma.

{\em Proof of proposition \ref{charf} }  Identities (\ref{chifqpi}) follow from lemma \ref{qbinom} for $t=q_i,~~m=1-a_{ij},~~ c=\varepsilon_{ij}^\pi a_{ij}$ since set of solutions (\ref{c2}) always contains $\pm (m-1)$.

Motivated by relations (\ref{fqpi}) we suggest new realizations of the quantum group $U_q^R ({\frak g})$, one for each Coxeter element $s_\pi$.
Let 
$U_q^\pi ({\frak g})$ be the associative algebra over $\Bbb C$ with generators 
$e_i , f_i , L_i^{\pm 1} i=1, \ldots l$ subject to the relations:

\begin{equation}\label{sqgr}
\begin{array}{l}
L_iL_j=L_jL_i ,~~ L_iL_i^{-1}=L_i^{-1}L_i=1 \\
 \\
L_ie_j L_i^{-1}=q_i^{ \delta _{i,j}}e_j,~~ 
L_if_j L_i^{-1}=q_i^{- \delta _{i,j}}f_j\\
 \\
e_i f_j -q^{ c^\pi _{ij}} f_j e_i = \delta _{i,j}{K_i -K_i^{-1} \over q_i -q_i^{-1}} , \\
 \\
K_i=\prod_{j=1}^lL_j^{a_{ji}} , \\
 \\
\sum_{r=0}^{1-a_{ij}}(-1)^r q^{r c_{ij}^\pi}
\left[ \begin{array}{c} 1-a_{ij} \\ r \end{array} \right]_{q_i} 
(e_i )^{1-a_{ij}-r}e_j (e_i)^r =0 ,~ i \neq j , \\
 \\
\sum_{r=0}^{1-a_{ij}}(-1)^r q^{r c_{ij}^\pi}
\left[ \begin{array}{c} 1-a_{ij} \\ r \end{array} \right]_{q_i} 
(f_i )^{1-a_{ij}-r}f_j (f_i)^r =0 ,~ i \neq j .
\end{array}
\end{equation}

It follows that the map $\tau_q^\pi :F_q^\pi \rightarrow U_q^\pi ({\frak g});~~~e_i\mapsto e_i$ is a {\em natural} embedding of $F_q^\pi $ into $U_q^\pi ({\frak g})$.  
From now on we identify $F_q^\pi $ with the subalgebra in $U_q^\pi ({\frak g}) $ generated by $e_i ,i=1, \ldots l$.

\begin{theorem} \label{newreal}
For every integer--valued solution $n_{ij}\in {\Bbb Z},~i,j=1,\ldots ,l$ of equations 

\begin{equation}\label{eqpi}
d_in_{ji}-d_jn_{ij}=c^\pi_{ij}
\end{equation}

there exists an algebra
isomorphism $\psi_{\{ n\}} : U_q^\pi ({\frak g}) \rightarrow 
U_q^R ({\frak g})$ defined  by formulas:

\begin{equation}
\begin{array}{l}
\psi_{\{ n\}}(e_i)=q_i^{-n_{ii}}\prod_{p=1}^lL_p^{n_{ip}}X_i^+ ,\\
 \\
\psi_{\{ n\}}(f_i)=\prod_{p=1}^lL_p^{-n_{ip}}X_i^- , \\
 \\
\psi_{\{ n\}}(L_i)=L_i .
\end{array}
\end{equation}

\end{theorem}

{\em Proof} is provided by direct verification of defining relations (\ref{sqgr}). The most nontrivial part is to verify  deformed q--Serre relations (\ref{fqpi}). The defining relations of $U_q^R ({\frak g})$  imply the following relations for $\psi_{\{ 
n\}}(e_i)$,

\begin{equation}\label{relconstr}
\sum_{k=0}^{1-a_{ij}}(-1)^k 
\left[ \begin{array}{c} 1-a_{ij} \\ k \end{array} \right]_{q_i}
q^{k(d_in_{ji}- {d_j}n_{ij})}\psi_{\{ n\}}(e_i)^k\psi_{\{ n\}}e_j\psi_{\{ n\}}(e_i)^{1-a_{ij}-k} =0 ,
\end{equation}
for any $i\neq j$.
Now using equation (\ref{eqpi}) we arrive to relations (\ref{fqpi}).

\begin{remark} 
The general solution of equation (\ref{eqpi}) 
is given by

\begin{equation}\label{eq3}
n_{ji}=\frac 12 (\varepsilon_{ij}a_{ij} + \frac{s_{ij}}{d_i}),.
\end{equation}
where $s_{ij}=s_{ji}$. In order to show that integer solutions,
$n_{ij} \in {\Bbb Z}$, exist, we choose $s_{ij}=b_{ij}$. Then,
$n_{ji}=\frac 12 (\varepsilon_{ij}+1)a_{ij}$, and this is an integer.

\end{remark}

We call the algebra $U_q^\pi ({\frak g})$ the Coxeter realization of the quantum group $U_q^R ({\frak g})$ corresponding to the Coxeter element $s_\pi$.
The subalgebra $F_q^\pi $ and the character $\chi_q^\pi$ are  quantum counterparts of $U({\frak n})$ and of the non--singular character $\chi_0$ , respectively.

\section{Non--singular characters and affine quantum groups}

In this section we suggest new realizations of affine quantum
groups labeled by Coxeter elements, similar to those described
in the previous section for finite--dimensional quantum groups.

Let $\widehat{\frak g}={\frak g}((z))\stackrel {\cdot}{+}{\Bbb C}$ 
be the nontwisted affine Lie algebra corresponding to $\frak g$ and 
let ${\frak n}((z))\subset \widehat{\frak g}$ be the loop algebra 
of the nilpotent Lie subalgebra ${\frak n} \subset {\frak g}$. 

Let $\chi$ be the character of $\frak n$ which takes value $1$ on
all root generators of $\frak n$. $\chi$ has a unique extension
to the character $\widehat \chi$ of  ${\frak n}((z))$,
such that $\widehat \chi$ vanishes on the complement
$z^{-1} {\frak n}[[z^{-1}]] + z {\frak n}[[z]]$
of $\frak n$ in  ${\frak n}((z))$.

It is our goal to define quantum counterparts of
the algebra $U({\frak n}((z)))\subset U(\widehat{\frak g})$ and 
of the character $\widehat \chi$. We start with the new Drinfeld 
realization of quantum affine algebras 
generalizing the loop realization of 
affine Lie algebras.

Let $U_q^R({\widehat{\frak g}})$ be an associative algebra generated by elements 
$X^\pm _{i,r} , r \in {\Bbb Z} ,~ H_{i,r} , 
r \in {\Bbb Z}\backslash \{ 0 \} ,~ L_i^{\pm 1} , i=1, \ldots l ~,~ 
q^{\pm \frac c2}$. Put 

\[
\begin{array}{l}
X^\pm _i (u)=\sum_{r \in {\Bbb Z}}X^\pm _{i,r}u^{-r} , \\
 \\
\Phi^\pm _i (u)=\sum_{r=0}^\infty \Phi^\pm _{i,\pm r}u^{\mp r}=
K_i^{\pm 1} exp\left( \pm (q_i-q_i^{-1})\sum_{s=1}^\infty H_{i,\pm s}u^{\mp s}\right) , \\
 \\
K_i=\prod_{j=1}^lL_j^{a_{ji}}.
\end{array}
\]

In terms of the generating series the defining relations are \cite{nr},\cite{kh-t}:

\begin{equation}\label{affqgr}
\begin{array}{l}
L_iL_j=L_jL_i ~,~ L_iL_i^{-1}=L_i^{-1}L_i=1 , \\
 \\
L_iX_j^\pm (u) L_i^{-1}=q_i^{\pm \delta _{i,j}}X_j^\pm (u) , \\
 \\
\Phi^\pm _i (u)\Phi^\pm _j (v)=\Phi^\pm _j (v)\Phi^\pm _i (u)~,~
L_i\Phi_j^\pm (u) L_i^{-1}=\Phi_j^\pm (u) , \\
\\
q^{\pm \frac c2} \mbox{ are central } , q^{ \frac c2}q^{- \frac c2}=1 , \\
\\
\Phi^+ _i (u)\Phi^- _j (v)=
{g_{ij}(\frac{vq^c}{u}) \over g_{ij}(\frac{vq^{-c}}{u})}\Phi^- _j (v)\Phi^+ _i (u), \\
\\ 
\Phi^- _i (u)X_j^\pm (v)\Phi^- _i (u)^{-1}=
g_{ij}(\frac{uq^{\mp \frac c2}}{v})^{\pm 1}X_j^\pm (v) , \\
 \\
\Phi^+ _i (u)X_j^\pm (v)\Phi^+ _i (u)^{-1}=
g_{ij}(\frac{vq^{\mp \frac c2}}{u})^{\mp 1}X_j^\pm (v) , \\
 \\
(u-vq^{\pm b_{ij}})X_i^\pm (u)X_j^\pm (v)=(q^{\pm b_{ij}}u-v)X_j^\pm (v)X_i^\pm (u) , \\
 \\
X_i^+(u)X_j^-(v) -X_j^-(v)X_i^+(u) = 
{\delta _{i,j} \over q_i - q_i^{-1}}
\left( \delta (\frac{uq^{-c}}{v})\Phi^+ _i (vq^{ \frac c2})-
\delta (\frac{uq^{c}}{v})\Phi^- _i (uq^{ \frac c2}) \right) , \\
 \\
\sum_{\pi \in S_{1-a_{ij}}}\sum_{k=0}^{1-a_{ij}}(-1)^k 
\left[ \begin{array}{c} 1-a_{ij} \\ k \end{array} \right]_{q_i}\times \\
X_i^\pm (z_{\pi (1)})\ldots X_i^\pm (z_{\pi (k)})X_j^\pm (w) 
X_i^\pm (z_{\pi (k+1)})\ldots X_i^\pm (z_{\pi (1-a_{ij})}) =0 , ~i \neq j , \\
 \\
\mbox{ where } g_{ij}(z)={1-q^{ b_{ij}}z \over 1-q^{ -b_{ij}}z}q^{ -b_{ij}} \in {\Bbb C}[[z]].
\end{array}
\end{equation}
The generators $X^\pm _{i,r} , H_{i,r}$ correspond to the elements $X_i^\pm z^r , H_i z^r$ of the affine Lie algebra $\widehat{\frak g}$ in the loop realization (here $X_i^\pm , H_i$ are the Chevalley generators of ${\frak g}$).

Let $I_n , n>0$ be the left ideal in $U_q^R({\widehat{\frak g}})$ 
generated by 
$X^\pm _{i,r} , i=1, \ldots l , r\geq n$ and by all polynomials in 
$H_{i,r} , r>0 ,~ L_i^{\pm 1}$ of degrees greater  or equal to $n$ 
$(deg(H_{i,r})=r ,~ deg(L_i^{\pm 1})=0)$. 
The algebra $\widehat U_q^R({\widehat{\frak g}})$ ,

\[
\widehat U_q^R({\widehat{\frak g}})=\lim_{\leftarrow}  U_q^R({\widehat{\frak g}})/I_n \mbox{ (inverse limit) }.
\]
is called the restricted completion of $U_q^R({\widehat{\frak g}})$.
We fix $k \in {\Bbb C}$ and  denote by 
$\widehat U_q^R({\widehat{\frak g}})_k$ the quotient of 
$\widehat U_q^R({\widehat{\frak g}})$ by
 the ideal generated by $q^{\pm \frac c2}-q^{\pm \frac k2}$. 
Sometimes it is convenient to use the weight--type generators $Y_{i,r}$,

\[
Y_{i,r}=(a^r)^{-1}_{ik}H_{k,r} , \, a^r_{ij}=\frac 1r [ra_{ij}]_{q_i} .
\]

In the spirit of theorem \ref{newreal}
we introduce, for any
any set of complex numbers $n_{ij}^{\pm r} , i,j=1, \ldots l , r\in {\Bbb N}$
and  integer parameters $n_{ij} , i,j=1, \ldots l $, generating
series of the form,
\begin{equation} \label{e} 
e_i^{\{ n\}}(u)=q_i^{-n_{ii}}{\Phi^0 _i}^{\{ n\}}{\Phi^- _i (u)}^{\{ n\}}X_i^+(u){\Phi^+ _i (u)}^{\{ n\}} ,
\end{equation}
where
\begin{equation}\label{phi}
\begin{array}{l}
{\Phi^\pm _i (u)}^{\{ n\}}=exp\left( \sum_{r=1}^\infty Y_{j,\pm r}\, n_{ij}^{\pm r} u^{\mp r} \right)  , 
n_{ij}^{\pm r} \in {\Bbb C} , \\
 \\
{\Phi^0 _i}^{\{ n\}}=\prod_{j=1}^lL_j^{n_{ji}} , n_{ij} \in {\Bbb Z}.
\end{array}
\end{equation}
The fourier coefficients of $e_i^{\{ n\}}(u)$ are 
elements of $\widehat U_q^R({\widehat{\frak g}})_k$.

\begin{proposition} \label{Fn}
The generating functions $e_i^{\{ n\}}(u)$ 
satisfy the following commutation relations,

\begin{equation}\label{erel}
\begin{array}{l}
(u-vq^{ b_{ij}})F_{ji}(\frac vu )e_i^{\{ n\}} (u)e_j^{\{ n\}} (v)=(q^{ b_{ij}}u-v)F_{ij}( \frac uv )e_j^{\{ n\}} (v)e_i^{\{ n\}} (u) , \\
 \\
\sum_{\pi \in S_{1-a_{ij}}}\sum_{k=0}^{1-a_{ij}}(-1)^k 
\left[ \begin{array}{c} 1-a_{ij} \\ k \end{array} \right]_{q_i}
\prod_{p<q}F_{ii}(\frac {z_{\pi (q)}}{z_{\pi (p)}} )\prod_{r=1}^k F_{ji}(\frac {w}{z_{\pi (r)}}) \times \\
\prod_{s=k+1}^{1-a_{ij}}F_{ij}(\frac {z_{\pi (s)}}{w})
e_i^{\{ n\}} (z_{\pi (1)})\ldots e_i^{\{ n\}} (z_{\pi (k)})e_j^{\{ n\}}(w) \times \\
e_i^{\{ n\}} (z_{\pi (k+1)})\ldots e_i^{\{ n\}} (z_{\pi (1-a_{ij})}) =0 , \mbox{ for } i \neq j , \\
\end{array}
\end{equation}
where
\begin{equation}
\begin{array}{l}
F_{ij}(z)=q_j^{n_{ij}}exp( \sum_{r=1}^\infty ( 
(n_{ij}^{-r}-n_{ji}^{r})q^{-\frac{kr}{2}}- 
n_{ik}^{-r}n_{jl}^{r}r(B^r)^{-1}_{kl}(q^{kr}-q^{-kr})) z^r), \\
\\
B_{ij}^r=q^{rb_{ij}}-q^{-rb_{ij}} .
\end{array}
\end{equation}
\end{proposition}

{\em Proof.}
Proposition is proved by direct substitution of the expression
for the generating series $e_i^{\{ n\} }$ and using the
defining relations of the algebra $\widehat{U}_q^R(\widehat{\frak g})_k$.

The commutation relations stated in proposition \ref{Fn}
are affine counterparts of the modified quantum Serre
relations (\ref{fqpi}). We denote by $\widehat F_q^{\{ n\}}$ 
the subalgebra in $\widehat U_q^R({\widehat{\frak g}})_k$ generated by
Fourier coefficients of $e_i^{\{ n\}}(u) , i=1, \ldots l$. 

Suppose that the map 
$\widehat \chi_q^{\{ n\}}:\widehat F_q^{\{ n\}}\rightarrow {\Bbb C}$ 
defined on the generators by equation 
$\widehat \chi_q^{\{ n\}}(e_i^{\{ n\}}(u))=1$ is a character of 
the algebra $\widehat F_q^{\{ n\}}$. Then, relations (\ref{erel})
 imply the following equations for the formal power
series $F_{ij}(z)$:

\begin{equation}\label{1constr}
(z-q^{ b_{ij}})F_{ji}(z^{-1})=(q^{ b_{ij}}z-1)F_{ij}( z) , a_{ij}\neq 0 ,
\end{equation}

\begin{equation}\label{2constr}
F_{ji}(z^{-1})=F_{ij}( z) , a_{ij}= 0 ,
\end{equation}

\begin{equation}\label{3constr}
\begin{array}{l}
\sum_{\pi \in S_{1-a_{ij}}}\sum_{k=0}^{1-a_{ij}}(-1)^k 
\left[ \begin{array}{c} 1-a_{ij} \\ k \end{array} \right]_{q_i}
\prod_{p<q}F_{ii}(\frac {z_{\pi (q)}}{z_{\pi (p)}} )\prod_{r=1}^k F_{ji}(\frac {w}{z_{\pi (r)}})\times \\
\prod_{s=k+1}^{1-a_{ij}}F_{ij}(\frac {z_{\pi (s)}}{w})=0 ,\,  a_{ij} \neq 0 ,~ i \neq j.
\end{array}
\end{equation}
Our aim is to solve this system of equations
with respect to parameters $n_{ij}, n_{ij}^r$. It is easy to
find solutions of
(\ref{1constr}) and (\ref{2constr}).

\begin{proposition}\label{p12}
Suppose that $d_jn_{ij}=d_in_{ji}$ for any $i$ and $j$ such that
$a_{ij}=0$. Then, the system of equations (\ref{1constr}) 
and (\ref{2constr}) has a unique solution
in Taylor series, ${\Bbb C}[[z]]$,  with constant terms 
$F_{ij}(0)=q_j^{n_{ij}}$.
This solution is given by formula,

\begin{equation}\label{F}
F_{ij}( z)={q_j^{n_{ij}} -zq_i^{n_{ji}} \over 1-zq^{b_{ij}}},~ a_{ij}\neq 0 ,
\end{equation}

\begin{equation}\label{F1}
F_{ij}(z)=q_j^{n_{ij}} ,~ a_{ij}=0 .
\end{equation}

\end{proposition}

Note that  parameters $n_{ij}$ in proposition \ref{p12} may take arbitrary complex values.
In applications to quantum groups $n_{ij}$
are integers.

{\em Proof.} Suppose that the Taylor series

\[
F_{ij}(z)=\sum_{n=0}^\infty c_{ij}^nz^n , 
\]
satisfy equations (\ref{1constr}). Then,
 the r.h.s. of (\ref{1constr}) are Taylor series. 
Therefore the l.h.s. of these equations must belong to ${\Bbb C}[[z]]$,
$(z-q^{ b_{ij}})F_{ji}(z^{-1}) \in {\Bbb C}[[z]]$. From the other hand
$(F_{ji}(z^{-1})-c_{ij}^0)(z-q^{ b_{ij}}) \in {\Bbb C}[[z^{-1}]]$. It follows 
that $(F_{ji}(z^{-1})-c_{ij}^0)(z-q^{ b_{ij}})=c_{ji} \in {\Bbb C}$, and
\begin{equation} \nonumber
F_{ji}(z^{-1})=c_{ij}^0+c_{ij}{z^{-1} \over 1-z^{-1}q^{b_{ij}}}.
\end{equation}
Substituting this anzatz into  (\ref{1constr}) we get the following relations for the 
coefficients $c_{ij}^0,~c_{ij}$:

$$
c_{ij}=-c_{ji}^0+q^{b_{ij}}c_{ij}^0 .
$$
Adding the condition $c_{ij}^0=q_j^{n_{ij}}$, one obtains (\ref{F}).

Equation (\ref{2constr}) implies that $F_{ij}(z)$ is a constant
for $i,j$ such that $a_{ij}=0$. Then, it is equal to the constant
term of its Taylor series, which gives (\ref{F1}).

Next, we show that, under some assumptions, equations 
(\ref{1constr}) and (\ref{2constr}) imply (\ref{3constr}).

\begin{proposition} \label{App}
Assume that $d_in_{ji}-d_jn_{ij}=c^\pi_{ij}$ for
some permutation $\pi \in S_l$. Then, any solution of
(\ref{F}) satisfies (\ref{3constr}).
\end{proposition}

{\em Proof.} We shall use theorem \ref{th1} and proposition \ref{th2} proved in 
Appendix. 

An important property of solution (\ref{F}) subject to the 
conditions of the proposition is that either 
$F_{ji}=q_i^{n_{ji}}$ or $F_{ij}=q_j^{n_{ij}}$. From  this fact it follows
that either the series in (\ref{constr1}) or the series in (\ref{constr2}) have
a common domain of convergence. Therefore either in 
(\ref{constr1}) or in (\ref{constr2}) the product of formal power series is  well defined. 
This allows us to apply theorem \ref{th1} or proposition \ref{th2}, respectively, to 
obtain identities (\ref{3constr}) for solution (\ref{F}).

Next, for some choice of  $n_{ij}$ as in proposition \ref{App}
we would like to choose complex parameters $n_{ij}^r, r \neq 0$
such that the following eqation is satisfied,
\begin{equation}\label{Kq}
(n_{ij}^{-r}-n_{ji}^{r})q^{-\frac{kr}{2}}- 
n_{ik}^{-r}n_{jl}^{r}\, r(B^r)^{-1}_{kl}(q^{kr}-q^{-kr})=
\frac 1r (q^{rb_{ij}}-q^{r(d_in_{ji}-d_jn_{ij})}) , r\in {\Bbb N} .
\end{equation}
Equation (\ref{Kq}) has many solutions. In particular, one
can choose $n_{ij}^r=0$ for all $r\le 0$. Then,
$$
n_{ji}^{r} = - q^{\frac{kr}{2}}
\frac 1r (q^{rb_{ij}}-q^{r(d_in_{ji}-d_jn_{ij})})
$$
for $r > 0$.

Finally, we conclude that
if coefficients $n_{ij}$ and $n_{ij}^{r}$ are solutions
of  equations (\ref{eqpi}) and (\ref{Kq}), the map
$\widehat \chi_q^{\{ n\}}$ is a character of the subalgebra $\widehat F_q^{\{ n\}}$.
Now we bring this result into focus by defining new realizations
of affine quantum algebras, similar to the new realizations
of finite--dimensional quantum algebras of the previous section.

First, we define quantum counterparts of $U({\frak n}((z)))$. 
Let $E_q$ be the free associative algebra generated by 
Fourier coefficients of the generating series $e_i (u) , i=1, \ldots l$.
Let $K_n , n>0$ be the left ideal in $E_q$ generated by $e_{i,r}, r\geq n$. 
Put $\widehat {E}_q=\lim_{\leftarrow}E_q/K_n$. Let 
$\widehat F_q^{\pi}$ be the quotient of $\widehat {E}_q$ by the 
two--sided ideal generated by the Fourier coefficients of 
the following generating series:

\begin{equation}
\begin{array}{l}
(u-vq^{ b_{ij}})F_{ji}(\frac vu )e_i (u)e_j (v)-(q^{ b_{ij}}u-v)F_{ij}( \frac uv )e_j (v)e_i (u) , \\
 \\
\sum_{\pi \in S_{1-a_{ij}}}\sum_{k=0}^{1-a_{ij}}(-1)^k 
\left[ \begin{array}{c} 1-a_{ij} \\ k \end{array} \right]_{q_i}
\prod_{p<q}F_{ii}(\frac {z_{\pi (q)}}{z_{\pi (p)}} )\prod_{r=1}^k F_{ji}(\frac {w}{z_{\pi (r)}}) \times \\
\prod_{s=k+1}^{1-a_{ij}}F_{ij}(\frac {z_{\pi (s)}}{w})
e_i (z_{\pi (1)})\ldots e_i (z_{\pi (k)})e_j(w) 
e_i (z_{\pi (k+1)})\ldots e_i (z_{\pi (1-a_{ij})})  , i \neq j ,
\end{array}
\end{equation}
where $F_{ij}(z)$ are given by (\ref{F}),(\ref{F1}) .

The algebra $\widehat F_q^{\pi}$ is an abstract version of subalgebras
$\widehat F_q^{\{ n\}}$. Note that the defining relations
of the algebra $\widehat F_q^{\pi}$ only depend on 
skew--symmetric combination (\ref{eqpi}) 
of the coefficients $n_{ij}$. Hence, similar to the
finite-dimensional case, one can associate 
$\widehat F_q^{\pi}$ to a Coxeter element of the Weyl group.

Next, we define new realizations of affine quantum groups.
Let $A_q$ be the free associative algebra generated by the Fourier coefficients of 
generating series 

\begin{equation}
\begin{array}{l}
e_i(u)=\sum_{r \in {\Bbb Z}}e_{i,r}u^{-r} ,\\
 \\
f_i(u)=\sum_{r \in {\Bbb Z}}f_{i,r}u^{-r} ,\\
 \\ 
K_i^\pm(u)=\sum_{r=0}^\infty K_{i,\pm r}^\pm u^{\mp r} ,\\
 \\ 
{K_i^\pm(u)}^{-1}=\sum_{r=0}^\infty {K_{i,\pm r}^\pm}^{-1}u^{\mp r} 
\end{array}
\end{equation}

and elements $L_i^{\pm 1} , i=1,\ldots ,l$.

Let $J_n , n>0$ be the left ideal in $A_q$ generated by $e_{i,r},f_{i,r}, r\geq n$ 
and by all polynomials in $K_{i,r}^+ ,{K_{i,r}^+}^{-1} ,L_i^{\pm 1} , r \geq 0$ of 
degrees greater than or equal to $n$ $(deg({K_{i,r}^+}^{\pm 1})=r , deg(L_i^{\pm 1})=0)$. 
Put $\widehat {A}_q=\lim_{\leftarrow}A_q/J_n$. 

For $F_{ij}(z)$ given by (\ref{F}),(\ref{F1}) we 
define the following formal power series: 
$$
\begin{array}{l}
M_{ij}(z)=g_{ij}(zq^{-k})^{-1}F_{ji}(zq^k)F_{ji}(zq^{-k})^{-1} , \\
 \\
G_{ij}(z)=M_{ij}(zq^{-k})M_{ij}(zq^k)^{-1} , \\
\\
F_{ij}^-(z)=F_{ij}(zq^{2k}).
\end{array}
$$

Let 
$\widehat U^\pi_{q,k}({\widehat{\frak g}})$ be the quotient of $\widehat {A}_q$ by the two--sided 
ideal generated by the Fourier coefficients of the following generating series:

$$
\begin{array}{l}
K_i^\pm(u)K_j^\pm(v)-K_j^\pm(v)K_i^\pm(u)~ , ~ K_i^\pm(u){K_i^\pm(u)}^{-1}-1~ ,~
{K_i^\pm(u)}^{-1}K_i^\pm(u)-1 ~, \\
 \\
L_iL_j-L_jL_i ~,~ L_iL_i^{-1}-1 ~,~ L_i^{-1}L_i-1 ~,\\
 \\
L_iK_j^\pm(v) L_i^{-1}-K_j^\pm(v) , \\
 \\
K_{i,0}^\pm - (\prod_{j=1}^lL_j^{a_{ji}})^{\pm 1} , \\
\\
K^+ _i (u)K^- _j (v)-G_{ij}(\frac vu )K^- _j (v)K^+ _i (u), 
\end{array}
$$

\begin{equation}\label{affpiqgr}
\begin{array}{l}
L_ie_j (u) L_i^{-1}-q_i^{ \delta _{i,j}}e_j (u) , \\
\\
L_if_j (u) L_i^{-1}-q_i^{- \delta _{i,j}}f_j (u) , \\
 \\
K^+ _i (u)e_j(v)-M_{ij}(\frac vu )e_j(v)K^+ _i (u), \\
\\
K^+ _i (u)f_j(v)-M_{ij}(\frac {vq^k}{u} )^{-1}f_j(v)K^+ _i (u), \\
 \\
K^- _i (u)e_j(v)-M_{ji}(\frac uv )^{-1}e_j(v)K^- _i (u),  \\
\\
K^- _i (u)f_j(v)-M_{ji}(\frac {uq^k}{v} )f_j(v)K^- _i (u),\\
\\
(u-vq^{ b_{ij}})F_{ji}(\frac vu )e_i (u)e_j (v)-(q^{ b_{ij}}u-v)F_{ij}( \frac uv )e_j (v)e_i (u) , \\
\\
(u-vq^{ -b_{ij}})F_{ji}^-(\frac vu )f_i (u)f_j (v)-(q^{ -b_{ij}}u-v)F_{ij}^-( \frac uv )f_j (v)f_i (u) , \\
\\
\sum_{\pi \in S_{1-a_{ij}}}\sum_{k=0}^{1-a_{ij}}(-1)^k 
\left[ \begin{array}{c} 1-a_{ij} \\ k \end{array} \right]_{q_i}
\prod_{p<q}F_{ii}(\frac {z_{\pi (q)}}{z_{\pi (p)}} )\prod_{r=1}^k F_{ji}(\frac {w}{z_{\pi (r)}})\times \\ 
\prod_{s=k+1}^{1-a_{ij}}F_{ij}(\frac {z_{\pi (s)}}{w}) 
e_i (z_{\pi (1)})\ldots e_i (z_{\pi (k)})e_j(w) 
e_i (z_{\pi (k+1)})\ldots e_i (z_{\pi (1-a_{ij})})  ,~ i \neq j , \\
 \\
\sum_{\pi \in S_{1-a_{ij}}}\sum_{k=0}^{1-a_{ij}}(-1)^k 
\left[ \begin{array}{c} 1-a_{ij} \\ k \end{array} \right]_{q_i}
\prod_{p<q}F_{ii}^-(\frac {z_{\pi (q)}}{z_{\pi (p)}} )\prod_{r=1}^k F_{ji}^-(\frac {w}{z_{\pi (r)}}) \times \\ 
\prod_{s=k+1}^{1-a_{ij}}F_{ij}^-(\frac {z_{\pi (s)}}{w}) 
f_i (z_{\pi (1)})\ldots f_i (z_{\pi (k)})f_j(w) 
f_i (z_{\pi (k+1)})\ldots f_i (z_{\pi (1-a_{ij})})  ,~ i \neq j , \\
\\
q_i^{n_{ji}}F_{ji}^{-1}(\frac {vq^k}{u})e_i(u)f_j(v) - 
q_i^{n_{ji}}F_{ij}^{-1}(\frac {uq^k}{v})f_j(v)e_i(u) - \\
-{\delta _{i,j} \over q_i - q_i^{-1}}
\left( \delta (\frac{uq^{-k}}{v})K^+ _i (v)-
\delta (\frac{uq^{k}}{v})K^- _i (v) \right) . 
\end{array}
\end{equation}
These relations only depend on skew--
symmetric combination (\ref{eqpi}) of the coefficients $n_{ij}$ . 
Thus, there is a one--to--one correspondence between
 Coxeter elements $s_\pi$ and the algebras 
$\widehat U^\pi_{q,k}({\widehat{\frak g}})$.

Of course, the algebra $\widehat F_q^{\pi}$ is a subalgebra
of $\widehat U^\pi_{q,k}({\widehat{\frak g}})$, with respect
to the natural embedding.

Finally, we show that $\widehat U^\pi_{q,k}({\widehat{\frak g}})$
is indeed a realization of $\widehat U_q^R({\widehat{\frak g}})_k$.

\begin{proposition}
For every integer--valued solution of equation (\ref{eqpi}) and every solution of 
system (\ref{Kq}) with $d_in_{ji}-d_jn_{ij}=c^\pi_{ij}$ there
exists an isomorphism of algebras $\widehat \psi_{\{ n\}} : \widehat U^\pi_{q,k}({\widehat{\frak g}}) 
\rightarrow \widehat U_q^R({\widehat{\frak g}})_k$ given by :

\begin{equation}
\begin{array}{l}
\widehat \psi_{\{ n\}}(e_i(u))=
q_i^{-n_{ii}}{\Phi^0 _i}^{\{ n\}}{\Phi^- _i (u)}^{\{ n\}}X_i^+(u)
{\Phi^+ _i (u)}^{\{ n\}} , \\
 \\
\widehat \psi_{\{ n\}}(f_i(u))=
{{\Phi^0 _i}^{\{ n\}}}^{-1}{{\Phi^- _i (uq^k)}^{\{ n\}}}^{-1}X_i^-(u)
{{\Phi^+ _i (uq^{-k})}^{\{ n\}}}^{-1} , \\
 \\
\widehat \psi_{\{ n\}}(K_i^\pm(u))=K_i^{\pm 1} exp( \sum_{s=1}^\infty 
\pm (q_i-q_i^{-1}) H_{i,\pm s}q^{-\frac {sk}{2}}u^{\mp s}- \\
\\
Y_{j,\pm s}n_{ij,\pm s}(q^{ks}-q^{-ks})u^{\mp s}) , \\
 \\
K_i=\prod_{j=1}^lL_j^{a_{ji}} , \\
 \\
\widehat \psi_{\{ n\}}(L_i)=L_i ,
\end{array}
\end{equation}

where ${\Phi^0 _i}^{\{ n\}} , {\Phi^- _i (u)}^{\{ n\}} , 
{\Phi^+ _i (u)}^{\{ n\}}$ are defined by (\ref{phi}).

\end{proposition}

{\em Proof.} The proof is by straightforward verification
of defining relations.

We shall identify $\widehat F_q^{\pi}$ with 
the subalgebra in $\widehat U^\pi_{q,k}({\widehat{\frak g}})$ generated by $e_{i,r} , 
r \in {\Bbb Z}$. Let $\widehat \chi_q^{\pi}$ be its canonical character ,
$\widehat \chi_q^{\pi}(e_i(u))=1$. For a fixed Coxeter element and for every 
solution of the corresponding equations (\ref{eqpi}) , (\ref{Kq}) the 
isomorphism $\widehat \psi_{\{ n\}}$ maps the subalgebra $\widehat F_q^{\pi}$ 
onto $\widehat F_q^{\{ n\}}$. The algebra $\widehat F_q^{\pi}$ and the 
character $\widehat \chi_q^{\pi}$ may be regarded as well--defined quantum 
counterparts of the algebra $U({\frak n}((z)))$ and the character 
$\widehat \chi$ , respectively. We shall call 
$\widehat U^\pi_{q,k}({\widehat{\frak g}})$  the Coxeter realization of 
the quantum group $\widehat U_q^R({\widehat{\frak g}})_k$ corresponding to 
$s_\pi$.

Finally observe that the map $\widehat \chi^{\pi}_\varphi : \widehat F_q^{\pi}
\rightarrow {\Bbb C}$ defined by 
$\widehat \chi^{\pi}_\varphi (e_i(u))=\varphi_i(u)¬,¬i=1,\ldots , l$, 
where $\varphi_i(u)\in {\Bbb C}((u))$ are arbitrary formal power series, is a character of 
the algebra $\widehat F_q^{\pi}$.

\section*{Appendix \\ A family of combinatorial identities}

In \cite{J} Jing proves the following identities for skew--symmetric polynomials.

\begin{theorem}(\cite{J})\label{l2}

For any $m\in {\Bbb Z} , m \leq 0$ the following identity holds:

\begin{equation}\label{trel}
\begin{array}{l}
\sum_{\pi \in S_{1-m}}(-1)^{l(\pi )}\sum_{k=0}^{1-m} 
\left[ \begin{array}{c} 1-m \\ k \end{array} \right]_{t}
\prod_{p<q}({z_{\pi (q)}}-t^2{z_{\pi (p)}} ) \times \\
\prod_{r=1}^k (1-t^m\frac {z_{\pi (r)}}{w})
\prod_{s=k+1}^{1-m}(\frac {z_{\pi (s)}}{w}-t^m)=0 .
\end{array}
\end{equation}

\end{theorem}
These identities were proved by looking at the representations
of current algebras, the negative integers $m$ arise as  
off-diagonal matrix entries $a_{ij}$ of Cartan matrixes.
Using Jing's identities (\ref{trel}) we will show that equations (\ref{1constr})
imply identities (\ref{3constr}). The latter are  generalizations
of Jing's identities.

Let $a_{ij} , i,j=1,\ldots , l$ be a generalized Cartan matrix : 
$a_{ii}=2$ , $a_{ij}$ are nonpositive integers for $i \neq j$, and 
$a_{ij}=0$ implies $a_{ji}=0$. Suppose 
also that $a_{ij}$ is symmetrizable , i. e. , there exist coprime positive 
integers $d_1,\ldots ,d_l$ such that the matrix $b_{ij}=a_{ij}d_i$ is 
symmetric. 

We shall make use of formal power series ( f.p.s. ) which are infinite in both 
directions. The space of such series is denoted by ${\Bbb C}((z))$. The product of 
two f.p.s. $f(z)=\sum_{n=-\infty}^{\infty}f_nz^n , g(z)=\sum_{n=-\infty}^{\infty}g_nz^n$ is 
said to exist if the coefficients of the series

\[
\sum_{p=-\infty}^{\infty}z^p\sum_{k+n=p}f_ng_k
\]
are well defined , i.e. the series
$\sum_{k+n=p}f_ng_k$ converges for every p.
Similarly , the product of three f.p.s $f(z)=\sum_{n=-\infty}^{\infty}f_nz^n$ , 
$g(z)=\sum_{n=-\infty}^{\infty}g_nz^n $ , $ h(z)=\sum_{n=-\infty}^{\infty}h_nz^n $ exists if the 
series $ \sum_{k+n+l=p}f_ng_kh_l$
converges for every p and its sum does not depend on the ordering of the terms.
Clearly , in this case the products $g(z)h(z) , g(z)f(z)$ and $f(z)h(z)$ are 
well--defined. For instance , if two or more formal power series
have a common domain of convergence their 
product is well--defined. We often use notation $\frac{1}{1-x}$ for the geometric series 
\[
\frac{1}{1-x}=\sum_{n=0}^\infty x^n
\]
viewed as a formal power series.

We show that equations (\ref{1constr}):

\begin{equation}\label{1}
(z-q^{ b_{ij}})F_{ji}(z^{-1})=(q^{ b_{ij}}z-1)F_{ij}( z) ,~ a_{ij}\neq 0 .
\end{equation} 
imply the  system of
identities (\ref{3constr}) for formal power series. Note that these identities 
hold true not only for the Taylor series solution of (\ref{1}).

\begin{theorem}\label{th1}
Let $F_{kl}(z),k,l=1,\ldots , l$ be a solution of equations (\ref{1}). 
Suppose that for some i and j, $i\neq j$, the following product is well-defined
as a formal power series,

\begin{equation}\label{constr1}
\prod_{p\neq q} {1 \over 1-q^{b_{ii}} \frac {z_q}{z_p}} \cdot  \prod_{s=1}^{1-a_{ij}}
{1 \over 1-q^{b_{ij}}\frac {z_s}{w}} \cdot P_{ij},
\end{equation}
where
\begin{equation} \label{2}
P_{ij} = 
\sum_{\pi \in S_{1-a_{ij}}}\sum_{k=0}^{1-a_{ij}}(-1)^k 
\left[ \begin{array}{c} 1-a_{ij} \\ k \end{array} \right]_{q_i}
\prod_{p<q}F_{ii}(\frac {z_{\pi (q)}}{z_{\pi (p)}} )
\prod_{r=1}^k F_{ji}(\frac {w}{z_{\pi (r)}})
\prod_{s=k+1}^{1-a_{ij}}F_{ij}(\frac {z_{\pi (s)}}{w}).
\end{equation}
Then, $P_{ij}=0$, as a formal power series.

\end{theorem}
The l.h.s of (\ref{2}) is symmetric with respect to permutations 
of the formal variables $z_1, \ldots , z_{1-a_{ij}}$. Thus, the last
theorem yields a family of combinatorial identities for symmetric functions.

Now we turn to the proof of the theorem.
First, we prove the following lemma.

\begin{lemma}\label{l1}
Let $F_{kl}(z),k,l=1,\ldots , l$ be a solution of system (\ref{1}). Then 

\begin{equation}\label{*}
\prod_{p\neq q}(1-q^{b_{ii}}\frac {z_q}{z_p})\prod_{s=1}^{1-a_{ij}}
(1-q^{b_{ij}}\frac {z_s}{w})P_{ij}=0 ,\mbox{ for } i \neq j .
\end{equation}

\end{lemma}

{\em Proof.} Let $\pi \in S_{1-a_{ij}}$. Consider the product:

\[
\prod_{s=1}^{1-a_{ij}}(1-q^{b_{ij}}\frac {z_s}{w})\prod_{r=1}^k F_{ji}(\frac {w}{z_{\pi (r)}})
\prod_{s=k+1}^{1-a_{ij}}F_{ij}(\frac {z_{\pi (s)}}{w}).
\]

The f.p.s. $\prod_{s=1}^{1-a_{ij}}(1-q^{b_{ij}}\frac {z_s}{w})$ is symmetric with 
respect to permutations of the formal variables $z_s$. Therefore

\[
\prod_{s=1}^{1-a_{ij}}(1-q^{b_{ij}}\frac {z_s}{w})\prod_{r=1}^k F_{ji}(\frac {w}{z_{\pi (r)}})
\prod_{s=k+1}^{1-a_{ij}}F_{ij}(\frac {z_{\pi (s)}}{w})=
\]

\[
\prod_{s=1}^{1-a_{ij}}(1-q^{b_{ij}}\frac {z_\pi (s)}{w})\prod_{r=1}^k F_{ji}(\frac {w}{z_{\pi (r)}})
\prod_{s=k+1}^{1-a_{ij}}F_{ij}(\frac {z_{\pi (s)}}{w}).
\]

Now using equations (\ref{1}) for $F_{ij}$ we obtain:

\[
\prod_{s=1}^{1-a_{ij}}(1-q^{b_{ij}}\frac {z_\pi (s)}{w})\prod_{r=1}^k F_{ji}(\frac {w}{z_{\pi (r)}})
\prod_{s=k+1}^{1-a_{ij}}F_{ij}(\frac {z_{\pi (s)}}{w})=
\]

\[
\prod_{s=1}^k(1-q^{b_{ij}}\frac {z_{\pi (s)}}{w})\prod_{s=k+1}^{1-a_{ij}}(q^{b_{ij}}-\frac {z_{\pi (s)}}{w})
\prod_{r=1}^{1-a_{ij}}F_{ji}(\frac {w}{z_{\pi (r)}})=
\]

\begin{equation}\label{3}
\prod_{s=1}^k(1-q^{b_{ij}}\frac {z_{\pi (s)}}{w})\prod_{s=k+1}^{1-a_{ij}}(q^{b_{ij}}-\frac {z_{\pi (s)}}{w})
\prod_{r=1}^{1-a_{ij}}F_{ji}(\frac {w}{z_{r}}),
\end{equation}

since $\prod_{r=1}^{1-a_{ij}}F_{ji}(\frac {w}{z_{\pi (r)}})$ is also a symmetric f.p.s..

Similarly,

\begin{equation}\label{4}
\begin{array}{l}
\prod_{p\neq q}(1-q^{b_{ii}}\frac {z_q}{z_p})\prod_{p<q}F_{ii}(\frac {z_{\pi (q)}}{z_{\pi (p)}} )= \\
\\
\prod_{p>q}\left( \frac 1{z_q}(1-q^{b_{ii}}\frac {z_q}{z_p})F_{ii}(\frac {z_q}{z_p} ) \right)
\prod_{p<q}(-1)^{l(\pi )}(z_{\pi (q)}-q^{b_{ii}}z_{\pi (p)}).
\end{array}
\end{equation}

Substituting (\ref{3}) and (\ref{4}) into (\ref{*}) we get :

\begin{equation}
\prod_{p\neq q}(1-q^{b_{ii}}\frac {z_q}{z_p})\prod_{s=1}^{1-a_{ij}}
(1-q^{b_{ij}}\frac {z_s}{w})P_{ij}=
\end{equation}

\begin{equation}\nonumber
\begin{array}{l} 
\prod_{r=1}^{1-a_{ij}}F_{ji}(\frac {w}{z_{r}})
\prod_{p>q}\left( \frac 1{z_q}(1-q^{b_{ii}}\frac {z_q}{z_p})F_{ii}(\frac {z_q}{z_p} ) \right) \times \\
\\
\sum_{\pi \in S_{1-a_{ij}}}\sum_{k=0}^{1-a_{ij}}(-1)^k 
\left[ \begin{array}{c} 1-a_{ij} \\ k \end{array} \right]_{q_i} \times \\
\\
\prod_{p<q}(-1)^{l(\pi )}(z_{\pi (q)}-q^{b_{ii}}z_{\pi (p)})\times \\
\\
\prod_{s=1}^k(1-q^{b_{ij}}\frac {z_{\pi (s)}}{w})
\prod_{s=k+1}^{1-a_{ij}}(q^{b_{ij}}-\frac {z_{\pi (s)}}{w}).
\end{array}
\end{equation}

Now lemma \ref{l1} follows immediately from
theorem \ref{l2} with $t = q_i,~ m=a_{ij}$ .

{\em Proof of the theorem.} The conditions of the theorem imply that,
as a formal power series, 

\begin{equation} \nonumber 
P_{ij} =  \prod_{p\neq q} {1 \over 1-q^{b_{ii}} \frac {z_q}{z_p}} \prod_{s=1}^{1-a_{ij}}
{1 \over 1-q^{b_{ij}}\frac {z_s}{w}}
 \prod_{p\neq q}(1-q^{b_{ii}}\frac {z_q}{z_p})\prod_{s=1}^{1-a_{ij}}
(1-q^{b_{ij}}\frac {z_s}{w})P_{ij}.
\end{equation}

From lemma  \ref{l1} it follows that

\[
\prod_{p\neq q}(1-q^{b_{ii}}\frac {z_q}{z_p})\prod_{s=1}^{1-a_{ij}}
(1-q^{b_{ij}}\frac {z_s}{w})P_{ij}=0.
\]
Therefore, $P_{ij}=0$. This concludes the proof.

One can formulate several versions of theorem \ref{th1}. For instance, the following is true.

\begin{proposition}\label{th2}
Let $F_{kl}(z),k,l=1,\ldots , l$ be a solution of system (\ref{1}). 
Suppose that for some i and j the product

\begin{equation}\label{constr2}
\prod_{p\neq q} {1 \over 1-q^{b_{ii}}\frac {z_q}{z_p}} \cdot \prod_1^{1-a_{ij}} 
\frac{\frac{w}{z_s}}{1-q^{b_{ij}}\frac{w}{z_s}} \cdot P_{ij}
\end{equation}
is well--defined as a formal power series. Then, $P_{ij}=0$.
\end{proposition}
Proof of the proposition is similar to that of theorem \ref{th1}.

Similar statements exist for $|q|>1$.


\begin{thebibliography}{99}

\bibitem{Bur} Bourbaki N. , Groupes et algebras de Lie, Chap. 4,5,6 , Paris, Hermann (1968).

\bibitem{ChP} Chari V. , Pressley A. , A guide to quantum groups , Cambridge Univ. Press (1994).

\bibitem{nr}  Drinfeld V.G., A new realization of Yangians and quantized
affine algebras, {\em Sov. Math. Dokl.} {\bf 36} (1988).


\bibitem{D} Drinfeld V.G. , Quantum groups, Proc. Int. Congr. Math. Berkley , California,1986, Amer. Math. Soc. , Providence (1987) , p.p. 718-820.

\bibitem{FF}  Feigin B., Frenkel E., Affine Lie algebras at the critical
level and Gelfand-Dikii algebras, {\em Int. J. Mod. Phys.} {\bf A7}, suppl.
A1 (1992), 197-215 ;

Quantization of the Drinfeld--Sokolov reduction , {\em Phys. Lett.} {\bf B 246} (1990) , 75--81.


\bibitem{GR} Gasper G. , Rahman M. , Basic hypergeometric series , Cambridge Univ. Press (1990).

\bibitem{Kac} Kac V. G. , Infinite dimensional Lie algebras , Cambridge Univ. Press (1990).

\bibitem{kh-t}  Khoroshkin S.M., Tolstoy V.N. , On Drinfeld's realization of
quantum affine algebras, J. Geom. Phys. {\bf11} (1993), 445-452.

\bibitem{K}  Kostant B. , The principal three--dimensional subgroup and the Betti numbers of a complex simple Lie group , {\em Amer. J. Math.} {\bf 81} (1959) , 973-1032.

\bibitem{K1}  Kostant B. , On Whittaker vectors and representation
theory, {\em Inventiones Math.} {\bf 48} (1978) , 101-184.

\bibitem{K2}  Kostant B. , The solution to a generalized Toda lattice and representation theory , {\em Adv. in Math.} {\bf 34} (1979) , 195-338.

\bibitem{J} N. Jing , {\em Quantum Kac--Moody algebras and vertex representations}, 
q-alg/9802036 .

\bibitem{M} I. G. Macdonald , {\em Symmetric functions and Hall polynomials} , 2nd 
edition , Claredon Press , Oxford , 1995.

\end{thebibliography}
\end{document}